\documentclass[aos,nameyear,dvips]{arximspdf}
\usepackage{graphics}
\doi{10.1214/09-AOS746}
\volume{38}
\issue{3}
\pubyear{2010}
\firstpage{1320}
\lastpage{1340}

\makeatletter

\newtheorem{lemma}{Lemma}
\newtheorem{theorem}{Theorem}
\newproclaim{Remark}{Remark}

\makeatother

\begin{document}
\begin{frontmatter}

\title{Asymptotic distribution of conical-hull estimators of directional edges}
\runtitle{Asymptotic distribution of conical-hull estimators}

\begin{aug}
\author[A]{\fnms{Byeong U.} \snm{Park}\thanksref{t1,t3}\ead[label=e1]{bupark@stats.snu.ac.kr}},
\author[B]{\fnms{Seok-Oh} \snm{Jeong}\corref{}\thanksref{t2,t3}\ead[label=e2]{seokohj@hufs.ac.kr}} and
\author[C]{\fnms{L\'{e}opold} \snm{Simar}\thanksref{t3}\ead[label=e3]{leopold.simar@uclouvain.be}}
\runauthor{B. U. Park, S.-O. Jeong and L. Simar}
\affiliation{Seoul National University, Hankuk University of Foreign Studies
and Universit\'{e}~catholique de Louvain}
\address[A]{B. U. Park\\
Department of Statistics\\
Seoul National University\\
South Korea\\
\printead{e1}}
\address[B]{S.-O. Jeong\\
Department of Statistics\\
Hankuk University of Foreign Studies\\
South Korea\\
\printead{e2}}
\address[C]{L. Simar\\
Institut de statistique\\
Universit\'{e} catholique de Louvain\\
Belgium\\
\printead{e3}}
\end{aug}

\pdfauthor{Byeong U. Park, Seok-Oh Jeong, Leopold Simar}

\thankstext{t1}{Supported by KOSEF funded by the Korea government (MEST 2009-0052815).}

\thankstext{t2}{Supported by a National Research Foundation of Korea
Grant funded by the Korean Government (2009-0067024).}

\thankstext{t3}{Supported by the ``Interuniversity Attraction Pole,''
Phase VI (No. P6/03) from the Belgian Science Policy.}

\received{\smonth{3} \syear{2009}}
\revised{\smonth{7} \syear{2009}}

%
\begin{abstract}
Nonparametric data envelopment analysis (DEA) estimators have been
widely applied in analysis of productive efficiency. Typically they are
defined in terms of convex-hulls of the observed combinations of
$\mathrm{inputs}\times\mathrm{outputs}$ in a sample of enterprises. The
shape of the convex-hull relies on a hypothesis on the shape of the
technology, defined as the boundary of the set of technically
attainable points in the $\mathrm{inputs}\times\mathrm{outputs}$ space.
So far, only the statistical properties of the smallest convex
polyhedron enveloping the data points has been considered which
corresponds to a situation where the technology presents variable
returns-to-scale (VRS). This paper analyzes the case where the most
common constant returns-to-scale (CRS) hypothesis is assumed. Here the
DEA is defined as the smallest conical-hull with vertex at the origin
enveloping the cloud of observed points. In this paper we determine the
asymptotic properties of this estimator, showing that the rate of
convergence is better than for the VRS estimator. We derive also its
asymptotic sampling distribution with a practical way to simulate it.
This allows to define a bias-corrected estimator and to build
confidence intervals for the frontier. We compare in a simulated
example the bias-corrected estimator with the original conical-hull
estimator and show its superiority in terms of median squared error.
\end{abstract}

%

\end{frontmatter}

\section{Introduction}

Consider a convex set $\Psi$ in $\mathbb{R}_+^{p+1}$ which takes the form
\[
\Psi= \{(\mathbf{x},y) \in\mathbb{R}_+^{p+1} \dvtx0 \le y \le g(\mathbf{x})\},
\]
where $g$ is a nonnegative convex function defined on $\mathbb{R}_+^p$
such that $g(a\mathbf{x}) = a g(\mathbf{x})$ for all $a>0$. Suppose
that we have a random sample $(\mathbf{X}_i,Y_i)$ drawn from a
distribution which is supported on $\Psi$. In this paper, we are
interested in estimating the ``boundary'' function $g$ from the random
sample. In particular, we study the asymptotic distribution of the
estimator
%
\begin{equation}\label{dea-est-g}
{\hat{g}}(\mathbf{x}) = \max\{y>0\dvtx(\mathbf{x},y)\in\widehat{\Psi} \},
\end{equation}
where $\widehat{\Psi}$ is the convex-hull of the rays $\mathbf{R}_i\equiv\{
(\gamma\mathbf{X}_i,
\gamma Y_i)\dvtx\gamma\ge0\}$ for all sample points $(\mathbf{X}_i,Y_i)$.

The problem arises in an area of econometrics where one is interested
in evaluating the performance of an enterprise in terms of technical
efficiency. In this context, $\mathbf{X}_i$ is the observed input
vectors of the $i$th enterprise, $Y_i$ is its observed productivity
and $\Psi$ is the production set of technically feasible pairs of input
and output. The property that $g(a\mathbf{x})=ag(\mathbf{x})$ for all
$a>0$, or, equivalently, $\Psi=a\Psi$ for all $a>0$, is called
``constant returns-to-scale'' (CRS), and the commonly used estimator of
$\Psi$ in this case is the CRS-version of the data envelopment analysis
(DEA) estimator defined by
\[
\widehat{\Psi}_0 =
\Biggl\{(\mathbf{x},y)\in\mathbb{R}_+^{p+1}\dvtx\mathbf{x}\ge\sum_{i=1}^n
\gamma_i \mathbf{X}_i, y\le\sum_{i=1}^n \gamma_i Y_i  \mbox{ for some }
\gamma_i\ge0,i=1,\ldots,n \Biggr\}.
\]
In fact, $\widehat{\Psi}_0$ given above is nothing else than the
smallest convex set containing all the rays $\mathbf{R}_i$ and the
hyperplane $\{(\mathbf{x},0) \dvtx\mathbf{x}\in\mathbb{R}^p \}$. To see
this, suppose that $(\mathbf{x},y)$ belongs to $\widehat{\Psi}_0$.
Then, there exist $\gamma_i \ge0$ such that $\mathbf{x}\ge\sum_{i=1}^n
\gamma_i \mathbf{X}_i$ and $y \le\sum_{i=1}^n \gamma_i Y_i$. For these
constants $\gamma_i$, define
\[
\gamma_i^* = \gamma_i \biggl(\frac{y}{\sum_{j=1}^n \gamma_j Y_j} \biggr)
\le\gamma_i
\]
for $1 \le i \le n$. Then $\sum_{i=1}^{n} \gamma_i^* Y_i = y$. Since
$\mathbf{x} \ge\sum_{i=1}^n \gamma_i \mathbf{X}_i \ge\sum_{i=1}^n
\gamma_i^* \mathbf{X}_i$, we have $\mathbf{x}^* \equiv\mathbf{x}-
\sum_{i=1}^n \gamma_i^* \mathbf{X}_i \ge\mathbf{0}$. This
shows\vspace*{1pt}
$(\mathbf{x}, y) = \sum_{i=1}^n (\gamma_i^* \mathbf{X}_i, \gamma_i^*
Y_i) + (\mathbf{x}^*, 0)$. The estimator ${\hat{g}}$ defined in
(\ref{dea-est-g}) and the one based on $\widehat{\Psi}_0$ are identical
with probability tending to one if the density of $(\mathbf{X}_i, Y_i)$
is bounded away from zero in a neighborhood of the boundary point
$(\mathbf{x},g(\mathbf{x}))$.

The problem that we describe in the first paragraph can be generalized
to the case of vector-valued $\mathbf{y}\in\mathbb{R}^q$. This is
particularly important in the specific problem that we mention in the
above paragraph where productivity is typically measured in several
variables. For this, we consider a conical-hull of a convex set $A$ in
$\mathbb{R}_+^{p+q}$ which is given by
\begin{eqnarray*}
& & \Psi\equiv\{(\mathbf{x},\mathbf{y}) \in\mathbb{R}_+^{p+q}
\dvtx\mbox{there exists a
constant } a>0
\mbox{ such that } (a\mathbf{x},a\mathbf{y}) \in A\} \cup
\{\mathbf{0}\}.
\end{eqnarray*}
The set $\Psi$ is convex and satisfies the CRS condition
%
\begin{equation}\label{crs}
a\Psi= \Psi\qquad\mbox{for all } a>0.
\end{equation}
We are interested in estimating the ``directional edge'' of $\Psi$ in
the $\mathbf{y}$-space, defined by
\[
\lambda(\mathbf{x},\mathbf{y})=\sup\{\lambda>0\dvtx(\mathbf{x},\lambda\mathbf{y})\in\Psi\}
\]
using a random sample from a density supported on $\Psi$. In the case
where $q=1$, the directional edge is linked directly to the boundary
function $g$ by the identity $g(\mathbf{x}) = y\lambda(\mathbf{x},y)$.
We consider the estimator
%
\begin{equation}\label{dea-est}
\hat{\lambda}(\mathbf{x},\mathbf{y}) =
\sup\{\lambda>0\dvtx(\mathbf{x},\lambda\mathbf{y})\in \widehat{\Psi}
\},
\end{equation}
where $\widehat{\Psi}$ is the convex-hull of the rays
$\mathbf{R}_i\equiv\{ (\gamma\mathbf{X}_i,
\gamma\mathbf{Y}_i)\dvtx\gamma\ge0\}$ for all sample points
$(\mathbf{X}_i,\mathbf{Y}_i)$.

To date, nonparametric data envelopment analysis (DEA) estimators have
been discussed or applied in more than 1800 articles published in more
than 400 journals [see Gattoufi, Oral and Reisman (\citeyear{GOR04})
for a comprehensive bibliography]. DEA estimators are used to estimate
various types of productive efficiency of firms in a wide variety of
industries as well as governmental agencies, national economies and
other decision-making units. The estimators employ linear programming
methods, similar to the one appearing in (\ref{dea-est}), along the
lines of Charnes, Cooper and Rhodes (\citeyear{CCR78}) who popularized
the basic ideas of Farrell (\citeyear{F57}).

Typically these DEA estimators are indeed defined in terms of
convex-hulls of the combinations of
$\mathrm{inputs}\times\mathrm{outputs}$ $(\mathbf{X}_i,\mathbf{Y} _i)$
in a sample of firms. The shape of the convex-hull relies on a
hypothesis on the shape of the technology defined as the boundary of
the set $\Psi$ of technically attainable points in the
$\mathrm{inputs}\times\mathrm{outputs}$ space. So far, only the
statistical properties of the smallest convex polyhedron enveloping the
data points has been considered which corresponds to a situation where
the technology presents variable returns-to-scale (VRS). Convergence
results for DEA--VRS have been derived by Korostelev, Simar and Tsybakov
(\citeyear{KST95}) in the case of univariate output and by Kneip, Park
and Simar (\citeyear{KPS98}) in the multivariate case. Asymptotic
distribution of the DEA--VRS estimators was obtained in the bivariate
case ($p=q=1$) by Gijbels et al. (\citeyear{Gijbelsetal99}), for
univariate output by Jeong and Park
(\citeyear{JP06}) and for the full multivariate case by Jeong (\citeyear {J04}) and  Kneip, Simar
and Wilson (\citeyear{KSW08}).

VRS is a flexible assumption, but in many situations the economist
assumes that the technology presents CRS: the first version of the DEA
estimator derived by Farrell (\citeyear{F57}) was for this situation.
Here the DEA estimator $\widehat\Psi$ is defined, as above, after (\ref
{dea-est}), as the smallest conical-hull with a vertex at the origin
enveloping the cloud of observed points. The properties of this
estimator have not been investigated, yet it was conjectured that one
would gain some efficiency in the estimation by imposing the
appropriate CRS structure to the estimator.

In this paper we determine the asymptotic properties of the DEA--CRS
estimator defined in (\ref{dea-est}), showing that the rate of
convergence is better than that of the VRS estimator. We derive also
its asymptotic sampling distribution with a practical way to simulate
it. This allows us to define a bias-corrected estimator and to build
confidence intervals for the frontier. We compare, in a simulated
example, the bias-corrected estimator with the original DEA--CRS
estimator and show its superiority in terms of median squared error.

\section{Rate of convergence}\label{sec2}

In this section we give the first theoretical result, the convergence rate
of the estimator $\hat{\lambda}$, as defined in (\ref{dea-est}),
in the general case of $p, q \ge1$. Before presenting the result, we
first give two lemmas which will be used in the proof of the first
theorem.
\begin{lemma}\label{lemma1}
For any $\alpha, \beta>0$, it holds that $\lambda (\alpha\mathbf{x},
\beta\mathbf{y}) = \frac{\alpha}{\beta}\lambda(\mathbf{x},\mathbf{y})$
whenever $(\alpha\mathbf{x}, \beta\mathbf{y}) \in\Psi$ and
$(\mathbf{x},\mathbf{y}) \in\Psi$. The same identity holds for
${\hat{\lambda}}$.
\end{lemma}
\begin{pf}
The lemma follows from the CRS property (\ref{crs})
since
\[
\sup\{\lambda>0\dvtx(\alpha\mathbf{x},\lambda\beta\mathbf{y})\in\Psi\}
= \sup\biggl\{\lambda>0\dvtx\biggl(\mathbf{x}, \frac{\lambda\beta}{\alpha}\mathbf{y}
\biggr)\in\Psi
\biggr\}.
\]
\upqed\end{pf}

The following lemma is also derived from the convexity of $\Psi$ and
$\widehat{\Psi}$.
\begin{lemma}\label{lemma2}
For all $r \in[0,1]$ and for all $(\mathbf{x}_1,\mathbf{y}_1),(\mathbf{x}_2,\mathbf{y}
_2)\in\Psi$,
\[
\lambda[r(\mathbf{x}_1,\mathbf{y}_1) + (1-r)(\mathbf{x}_2,\mathbf{y}_2) ] \ge
r \lambda(\mathbf{x}_1,\mathbf{y}_1)+(1-r)\lambda(\mathbf{x}_2,\mathbf{y}_2).
\]
The same inequality holds for ${\hat{\lambda}}$.
\end{lemma}

Our first theorem on the rate of convergence relies on the following
assumptions. In what follows, we fix the point in $\Psi$ where we want
to estimate $\lambda$, and denote it by $(\mathbf{x}_0, \mathbf{y}_0)$.
Throughout the paper, we assume that $(\mathbf{X}_i, \mathbf{Y}_i)$ are
independent and identically distributed with a density $f$ supported on
$\Psi\subset \mathbb{R}_+^p\times\mathbb{R}_+^q$ and that
$(\mathbf{x}_0,\mathbf{y}_0)$ is in the interior of $\Psi$.
\begin{enumerate}[(A1)]
\item[(A1)] $\lambda(\mathbf{x},\mathbf{y})$ is twice partially continuously
differentiable in a
neighborhood of $(\mathbf{x}_0,\mathbf{y}_0)$.
\item[(A2)] The density $f$ of $(\mathbf{X},\mathbf{Y})$ on
$\{(\mathbf{x},\mathbf{y})\in\Psi\dvtx\|(\mathbf{x},\mathbf{y}) -
(\mathbf{x}_0, \lambda(\mathbf{x}_0,\mathbf{y}_0) \mathbf{y} _0)\|
\le\varepsilon\}$ for some $\varepsilon>0$ is bounded away from zero.
\end{enumerate}
\begin{theorem}\label{theorem1}
Under the assumptions \textup{(A1)} and \textup{(A2)}, it follows that $
\hat{\lambda}(\mathbf{x}_0,\mathbf{y}_0)-\lambda(\mathbf{x}_0,\mathbf{y}_0)
= O_p(n^{-2/(p+q)})$.
\end{theorem}
\begin{pf}
We apply the technique of Kneip, Park and Simar (\citeyear{KPS98}). Put
$B_p(\mathbf{t},r)=\{\mathbf{x}\in\mathbb{R}_+^p\dvtx\|\mathbf{x}-\mathbf{t}\|\le
r\}$ and consider the balls near $\mathbf{x}_0\dvtx
C_r=B_p(\mathbf{x}_0^{(r)},\break h/2)$, $r=1,\ldots, 2p$ where
$\mathbf{x}_0^{(2j-1)}=\mathbf{x}_0-h\mathbf{e}_j$,
$\mathbf{x}_0^{(2j)}=\mathbf{x}_0+h \mathbf{e}_j$, $\mathbf{e}_j$ is
the unit $p$-vector with the $j$th element equal to 1 for
$j=1,2,\ldots,p$. Similarly, define $D_s=B_q(\mathbf{y}_0^{(s)},h/2)$
for $s=1,\ldots,2q$. Take $h$ small enough so that $C_r \times D_s
\subset \Psi$ for all $r=1,\ldots,2p$ and $s=1,\ldots,2q$. For
$r=1,\ldots,2p$, consider the conical hull of $C_r$,
\[
\mathcal{C}_r= \{\mathbf{x}\in\mathbb{R}_+^p \dvtx\exists a >0 \mbox{
such that } a \mathbf{x} \in C_r \}.
\]
Similarly, define $\mathcal{D}_s$. Define
\[
(\mathbf{U}_r,\mathbf{V}_s)=\mathop{\arg\min}_{(\mathbf{X}_i,\mathbf{Y}_i)
\in\mathcal{C}_r \times\mathcal{D}_s}
\lambda(\mathbf{X}_i,\mathbf{Y}_i).
\]
Since the number of points in $\mathcal{X}_n$ falling into $\Psi\cap
[\mathcal{C}_r\times\mathcal{D}_s]$ is proportional to $nh^{p+q-2}$, we have by
assumption (A2),
%
\begin{equation}\label{r}\quad
\lambda(\mathbf{U}_r,\mathbf{V}_s)=1 + O_p (n^{-1}h^{-p-q+2} ),\qquad r=1,\ldots
,2p, s=1,\ldots,2q.
\end{equation}

Let $\mathbf{U}_r^*= \alpha_r \mathbf{U}_r$ and $\mathbf{V}_s^*=\beta_s
V_s$ for $r=1,\ldots,2p$ and $s=1,\ldots,2q$ where $\alpha_r$ and
$\beta_s$ are positive constants such that $\mathbf{U}_r^* \in C_r$ and
$\mathbf{V}_s^* \in D_s$. Then from Lem\-ma~\ref{lemma1}, (\ref{r}) and
the fact that $\lambda, {\hat{\lambda}}\ge1$, it holds that for
$r=1,\ldots,2p$ and $s=1,\ldots,2q$,
\[
\frac{{\hat{\lambda}}(\mathbf{U}_r^*,\mathbf{V}_s^*)}{\lambda(\mathbf{U}_r^*,\mathbf{V}_s^*)}
=\frac{{\hat{\lambda}}(\mathbf{U}_r,\mathbf{V}_s)}{\lambda(\mathbf{U}_r,\mathbf{V}_s)}
\ge\frac{1}{\lambda(\mathbf{U}_r,\mathbf{V}_s)} = 1+O_p
(n^{-1}h^{-p-q+2} ),
\]
which implies that ${\hat{\lambda}}(\mathbf{U}_r^*,\mathbf{V}_s^*) \ge
\lambda(\mathbf{U}_r^*,\mathbf{V}_s^*)+O_p (n^{-1}h^{-p-q+2} )$. Since
$C_r$ and $D_s$ are balls surrounding the point $(\mathbf{x}_0,
\mathbf{y}_0)$, there exist scalars $w_r \ge0$ and $\omega_s \ge0$ such
that $\sum_{r=1}^{2p} w_r=1$, $\sum_{s=1}^{2q}\omega_s=1$,
$\mathbf{x}_0=\sum_{r=1}^{2p} w_r \mathbf{U}_r^*$ and
$\mathbf{y}_0=\sum_{s=1}^{2q} \omega_s \mathbf{V}_s^*$. Thus, from the
assumption (A1) we have
\[
\sum_{r=1}^{2p}\sum_{s=1}^{2q} w_r \omega_s
\lambda(\mathbf{U}_r^*,\mathbf{V}_s^*)=
\lambda(\mathbf{x}_0,\mathbf{y}_0) + O_p(h^2)
\]
for all $r$ and $s$. This, with Lemma \ref{lemma2} and the fact that
$\lambda\ge {\hat{\lambda}}$, shows that
\begin{eqnarray*}
\lambda(\mathbf{x}_0,\mathbf{y}_0) &\ge&
{\hat{\lambda}}(\mathbf{x}_0,\mathbf{y}_0)
\ge\sum_{r=1}^{2p}\sum_{s=1}^{2q} w_r \omega_s {\hat{\lambda}}(\mathbf{U}_r^*,\mathbf{V}_s^*)\\
&\ge& \sum_{r=1}^{2p}\sum_{s=1}^{2q} w_r \omega_s
\lambda(\mathbf{U}_r^*,\mathbf{V}
_s^*)+ O_p(n^{-1}h^{-p-q+2})\\
&=& \lambda(\mathbf{x}_0,\mathbf{y}_0) + O_p(h^2) + O_p (n^{-1}h^{-p-q+2} ).
\end{eqnarray*}
Taking $h \sim n^{-1/(p+q)}$ completes the proof of the theorem.
\end{pf}
\begin{Remark}\label{Remark1}
In the case where $\Psi$ is a convex set in $\mathbb{R}^{p+q}$
without having the CRS property (\ref{crs}), the DEA (data envelopment
analysis) estimator defined as in (\ref{dea-est}) with $\widehat{\Psi}$
replaced by the convex-hull of $(\mathbf{X}_i,\mathbf{Y}_i)$ is commonly used. In
this case, the DEA
estimator of $\lambda(\mathbf{x}_0,\mathbf{y}_0)$ is known to have $n^{-2/(p+q+1)}$
rate of convergence which is slightly worse than $n^{-2/(p+q)}$ [see
Kneip, Park and Simar (\citeyear{KPS98})]. The CRS property reduces the
``effective'' dimension by one.
\end{Remark}

\section{Asymptotic distribution}\label{sec3}

In this section we derive a representation for the asymptotic
distribution of the estimator ${\hat{\lambda}}$ defined in (\ref{dea-est}). This
representation enables one to simulate the asymptotic distribution so
that one can correct the bias of the estimator to get an improved
version of ${\hat{\lambda}}$. We work with the case where $q=1$ first and then
move to the general case where $q>1$. The result for the case $q=1$ is
essential for the generalization to $q>1$.

\subsection{The case where $q=1$}\label{sec31}

We consider the set
\[
\Psi= \{(\mathbf{x},y) \in A_c \times\mathbb{R}_+ \dvtx0 \le y \le g(\mathbf{x})\},
\]
where $g$ is a nonnegative convex function defined on a conical-hull
$A_c$ of a convex set $A \subset\mathbb{R}_+^p$ such that
%
\begin{equation}\label{crs-1}
g(a \mathbf{x}) = a g(\mathbf{x}) \qquad\mbox{for all } a >0,
\end{equation}
and that, for all $\mathbf{x}_1,\mathbf{x}_2\in A_c$ with $\mathbf{x}_1 \neq a\mathbf{x}_2$ for any $a>0$,
%
\begin{equation}\label{str_conv}
g\bigl(\alpha\mathbf{x}_1+(1-\alpha)\mathbf{x}_2\bigr) > \alpha g(\mathbf{x}_1)+(1-\alpha)g(\mathbf{x}_2)
\end{equation}
for all $\alpha\in(0,1)$. In this case, $\lambda(\mathbf{x}_0,y_0) = g(\mathbf{x}
_0)/y_0$ so that the problem of estimating $\lambda(\mathbf{x}_0,y_0)$ reduces
to that of estimating the function $g$ at $\mathbf{x}_0$. The estimator of
$g(\mathbf{x}_0)$ that corresponds to ${\hat{\lambda}}(\mathbf{x}_0,y_0)$ defined in (\ref
{dea-est}) is given by
%
%
\begin{equation}\label{g-est}
{\hat{g}}(\mathbf{x}_0) = y_0 {\hat{\lambda}}(\mathbf{x}_0,y_0) = \sup\{y \dvtx(\mathbf{x}_0,y) \in\widehat
{\Psi}\}.
\end{equation}
We note that the CRS condition (\ref{crs-1}) is satisfied, not only by linear
functions of the form $g(\mathbf{x}) = \mathbf{c}^\top\mathbf{x}$, but also by those
functions $g(\mathbf{x})=c (x_1^r + \cdots+ x_p^r)^{1/r}$ for all positive
numbers $c$
and positive integers $r$.

Define $\mathbf{S}_i$ by $\mathbf{S}_i^\top= (\mathbf{X}_i^\top, Y_i)$.
Below we describe a canonical transformation $T$ on $\Psi$ such that
the transformed data $T(\mathbf{S}_i)$ behave, asymptotically, as an
i.i.d. sample from a uniform distribution on a region that can be
represented by a simple $(p-1)$-dimensional quadratic function in the
transformed space. The reduction of the dimension, by one, for the
boundary function is due to the CRS property (\ref{crs-1}). This is
consistent with the dimension reduction as we noted in Remark
\ref{Remark1} in the previous section.

The key element in the derivation of the asymptotic distribution of
${\hat{g}} (\mathbf{x}_0)$ is to project the data $\mathbf{S}_i$ onto a
hyperplane which is perpendicular to the vector $\mathbf{x}_0$ and
passes through $\mathbf{x}_0$. The projected points lie under the locus
of the function $g$ on the hyperplane, and the estimator
${\hat{g}}(\mathbf{x}_0)$ equals the maximal $y$ such that
$(\mathbf{x}_0, y)$ belongs to the convex-hull of the projected points.
The asymptotic distribution of the estimator ${\hat{g}}(\mathbf{x}_0)$
is then obtained by analyzing the statistical properties of the
convex-hull of the projected points.

Let $Q$ be a $p \times(p-1)$ matrix whose columns constitute an
orthonormal basis for $\mathbf{x}_0^\perp$, the subspace of
$\mathbb{R}^{p}$ that is perpendicular to the vector $\mathbf{x}_0$.
Think of the transformation
\[
T_1
\dvtx\mathbf{x}\mapsto\biggl(\frac{\mathbf{x}_0^\top\mathbf{x}}{\|\mathbf{x}_0\|},
\mathbf{x}^\top Q \biggr)^\top.
\]
This transformation maps $\mathbf{x}$ to a vector which corresponds to
$\mathbf{x}$ in the new coordinate system where the axes are
$\mathbf{x}_0$ and the columns of $Q$. The first component of
$T_1(\mathbf{x})$ is nothing other than the projection of $\mathbf{x}$
onto the space spanned by $\mathbf{x}_0$, and the vector of the rest
components is its orthogonal complement in $\mathbb{R}^p$. Thus, the
inverse transform $T_1^{-1}$ is given by
\[
T_1^{-1} \dvtx\mathbf{z}\mapsto z_1
\biggl(\frac{\mathbf{x}_0}{\|\mathbf{x}_0\|} \biggr) + Q\mathbf{z}_2,
\]
where $\mathbf{z}^\top= (z_1, \mathbf{z}_2^\top)$.

It would be more convenient to use a transformation that takes $\mathbf{x}_0$
to the origin in the new coordinate system. This can be done by the
following transformation:
\[
T_2
\dvtx\mathbf{x}\mapsto\biggl[\frac{\mathbf{x}_0^\top(\mathbf{x}-\mathbf{x}_0)}{\|\mathbf{x}_0\|},
\biggl(\frac{\|\mathbf{x}_0\|^2}{\mathbf{x}_0^\top\mathbf{x}}
\biggr)\mathbf{x}^\top Q \biggr]^\top.
\]
Scaling by the factor $\|\mathbf{x}_0\|^2/\mathbf{x}_0^\top\mathbf{x}$
is introduced to factor out a common scalar for the inverse map of
$T_2$. In fact, $\|\mathbf{x}_0\|^2/\mathbf{x}_0^\top\mathbf{x}$ equals
the scalar $c$ such that the projection of $c\mathbf{x}$ onto the
linear span of $\mathbf{x}_0$ equals $\mathbf{x}_0$ itself. Thus
\[
\frac{\|\mathbf{x}_0\|^2}{\mathbf{x}_0^\top\mathbf{x}} \mathbf{x}= \mathbf{x}_0 + Q \biggl(Q^\top
\frac{\|\mathbf{x}_0\|^2}{\mathbf{x}_0^\top\mathbf{x}} \mathbf{x}\biggr)
\]
so that the inverse transform of $T_2$ is given by
\[
T_2^{-1} \dvtx\mathbf{z}\mapsto\biggl(\frac{z_1+\|\mathbf{x}_0\|}{\|\mathbf{x}_0\|} \biggr)
(\mathbf{x}_0 + Q\mathbf{z}_2 ).
\]
Note that $\mathbf{x}_0^\top\mathbf{x}>0$ if $\mathbf{x}\neq\mathbf{0}$
since then $\mathbf{x}_0, \mathbf{x}> \mathbf{0}$. It is easy to see
that $T_2(\mathbf{x}_0) = \mathbf{0}$.

Define a $(p-1)$-dimensional function $g^*$ by $g^*(\mathbf{z}_2) =
g(\mathbf{x}_0 + Q\mathbf{z}_2)$. For a function $\psi$, let
$\dot{\psi}$ and $\ddot{\psi}$ denote, respectively, the gradient
vector and the Hessian matrix of $\psi$. Since, for any
$\mathbf{u}\in\mathbb{R}^{p-1}$,
\[
\mathbf{u}^\top\ddot{g}^*(\mathbf{z}_2) \mathbf{u}=
(Q\mathbf{u})^\top\ddot{g}(\mathbf{x}_0 + Q\mathbf{z}_2) (Q\mathbf{u})
\]
and also $(Q\mathbf{u})^\top(Q\mathbf{u}) = \mathbf{u}^\top\mathbf{u}$,
it can be seen that $g^*$ is convex if $g$ is convex. In particular,
(\ref{str_conv}) implies the strict convexity of $g^*$. Note that $g^*$
does not have the CRS property (\ref{crs-1}), however.

Next, we introduce a further transformation on the new coordinate
system $(\mathbf{z},y)$. This transformation maps the equation
$y=g^*(\mathbf{z}_2)$ to a perfect quadratic equation in the further
transformed space. Since $g^*$ is strictly convex,
$-\ddot{g}^*(\mathbf{0})/2 = Q^\top (-\ddot{g}(\mathbf{x}_0)/2)Q$ is
positive definite and symmetric. Thus, there exist an orthogonal matrix
$P$ and a diagonal matrix $\Lambda$ such that
$-\ddot{g}^*(\mathbf{0})/2 = P \Lambda P^\top$. The columns of $P$ are
the orthonormal eigenvectors, and the diagonal elements of $\Lambda$
are the eigenvalues of the matrix $-\ddot{g}^*(\mathbf{0})/2$. Let
$T_3$ be a transformation that maps $\mathbb{R}^{p}$ to
$\mathbb{R}^{p}$ defined by
%
%
\begin{equation}\label{T3}
T_3 \dvtx\mathbf{z}\mapsto\bigl(z_1, n^{1/(p+1)} \mathbf{z}_2^\top P
\Lambda^{1/2} \bigr)^\top.
\end{equation}
Note that this transformation does not change $z_1$, the first
component of $\mathbf{z}$. Also, define a map $T_4\dvtx\mathbb{R}^{p}
\times\mathbb{R}\rightarrow\mathbb{R}$ by
%
%
\begin{equation}\label{T4}
T_4 \dvtx(\mathbf{z},y) \mapsto n^{2/(p+1)} \biggl[y
\biggl(\frac{\|\mathbf{x}_0\|}{z_1+\|\mathbf{x}_0\|} \biggr) -
g^*(\mathbf{0}) - \dot{g}^*(\mathbf{0})^\top\mathbf{z}_2 \biggr].
\end{equation}

The transformation we apply to the data $(\mathbf{X}_i,Y_i)$ is now defined by
\[
T \dvtx(\mathbf{x}, y) \mapsto\bigl(T_3 \circ T_2 (\mathbf{x}), T_4
(T_2 (\mathbf{x}),y ) \bigr).
\]
We explain how the equation $y=g(\mathbf{x})$ can be approximated,
locally at $(\mathbf{x}_0, y_0)$, by a $(p-1)$-dimensional quadratic
function in the new coordinate system transformed by $T$. Let
$(\mathbf{v},w) \in\mathbb{R}^p \times\mathbb{R}$ represent the new
coordinate system obtained by the transformation $T$. Write
$\mathbf{v}^\top=(v_1, \mathbf{v}_2^\top)$ with $\mathbf{v}_2$ being a
$(p-1)$-dimensional vector. Then, the inverse transform of $T$ maps
$\mathbf{v}$ and $w$, respectively, to
\begin{eqnarray*}
\mathbf{x}&=& \biggl(\frac{v_1 + \|\mathbf{x}_0\|}{\|\mathbf{x}_0\|} \biggr)
\bigl[\mathbf{x}_0 + n^{-1/(p+1)}Q P \Lambda^{-1/2} \mathbf{v}_2 \bigr],\\
y&=& \biggl(\frac{v_1 + \|\mathbf{x}_0\|}{\|\mathbf{x}_0\|} \biggr) \bigl[g^*(\mathbf{0}) +
n^{-1/(p+1)}\dot{g}^*(\mathbf{0})^\top P \Lambda^{-1/2}\mathbf{v}_2 + n^{-2/(p+1)}
w \bigr].
\end{eqnarray*}
Thus, for arbitrary compact sets $C_1 \subset\mathbb{R}^{p-1}$ and $C_2
\subset
\mathbb{R}$, we obtain using the CRS property (\ref{crs-1}) that, uniformly for
$v_1 \in\mathbb{R}_+$, $\mathbf{v}_2 \in C_1$ and $w \in C_2$,
\begin{eqnarray*}
&& y = g(\mathbf{x})\\
&&\quad\leftrightarrow\quad g^*(\mathbf{0}) + n^{-1/(p+1)}\dot{g}^*(\mathbf{0})^\top P
\Lambda^{-1/2} \mathbf{v}_2 + n^{-2/(p+1)} w\\
&&\hspace*{35.7pt}\qquad = g^* \bigl(n^{-1/(p+1)}P \Lambda^{-1/2} \mathbf{v}_2 \bigr)\\
&&\quad\leftrightarrow\quad w=-\mathbf{v}_2^\top\mathbf{v}_2 + o(1)
\end{eqnarray*}
as $n$ tends to infinity, provided that $\ddot{g}^*$ is continuous at
$\mathbf{0}$.

Now we give a representation of the limit distribution of ${\hat{g}}$ as given
in (\ref{g-est}). Define
%
%
\begin{eqnarray}\label{theta}
\theta&=& \|\mathbf{x}_0\| \int_0^\infty u^p f(u\mathbf{x}_0, ug(\mathbf{x}_0)) \,du,\\
\label{kappa}
\kappa&=& \theta\det(\Lambda)^{-1/2}.
\end{eqnarray}
Define a set $R_n(\kappa) \subset\mathbb{R}^{p}$ of points $(\mathbf{v}_2, w)$ such
that
\begin{eqnarray*}
\mathbf{v}_2 &\in& \bigl[ -\tfrac{1}{2}\kappa^{-1/(p+1)} n^{1/(p+1)},
\tfrac{1}{2}\kappa^{-1/(p+1)} n^{1/(p+1)} \bigr]^{p-1},\\
w &\in& \bigl[ -\mathbf{v}_2^\top\mathbf{v}_2 - \kappa^{-2/(p+1)}n^{2/(p+1)},
-\mathbf{v}_2^\top\mathbf{v}_2 \bigr].
\end{eqnarray*}
The volume of this set in $\mathbb{R}^p$ equals $n \kappa^{-1}$. Let $(\mathbf{V}_{2i},
W_i)$ be a random sample from the uniform distribution on $R_n(\kappa)$.
This random sample can be generated once we know $\kappa$. Let
$Z_n(\cdot)$ be defined as ${\hat{g}}$ in (\ref{g-est}) with $\widehat{\Psi}$
being replaced by the convex-hull of $(\mathbf{V}_{2i}, W_i)$; that said,
%
%
\begin{equation}\label{limit}\qquad
Z_n(\mathbf{v}_2)
= \sup\Biggl\{\sum_{i=1}^n \gamma_i W_i \dvtx\mathbf{v}_2 = \sum_{i=1}^n \gamma_i
\mathbf{V}_{2i}, \sum_{i=1}^n \gamma_i =1, \gamma_i\ge0, i=1,\ldots,n\Biggr\}.
\end{equation}
For a small $\varepsilon>0$, define a set on $\mathbb{R}_+^{p+1}$ by
%
%
\begin{eqnarray} \label{neigh}
&& H_\varepsilon(\mathbf{x}_0)
= \bigl\{ \bigl(u(\mathbf{x}_0+Q\mathbf{z}_2),
u\bigl(g(\mathbf{x}_0+Q\mathbf{z}_2)-y\bigr) \bigr)\dvtx u \ge0,\nonumber\\[-8pt]\\[-8pt]
&&\hspace*{154pt} \|\mathbf{z}_2\| \le\varepsilon, 0 \le y
\le\varepsilon\bigr\}.\nonumber
\end{eqnarray}

In the theorem below and those that follow, we will measure the
distance between two distributions by the following modification of the
Mallows distance:
\[
d(\mu_1, \mu_2) = \inf_{Z_1, Z_2} \{ E(Z_1-Z_2)^2 \wedge1 \dvtx\mathcal{L}
(Z_1)=\mu_1, \mathcal{L}(Z_2) = \mu_2 \}.
\]
Convergence in this metric is equivalent to weak convergence.
\begin{theorem}\label{theorem2}
Assume \textup{(A1)} and \textup{(A2)}. In addition, assume that $-\ddot{g}^*$
is positive definite and continuous at $\mathbf{0}$ and that the density $f$
of $(\mathbf{X}, Y)$ is uniformly continuous on $H_\varepsilon(\mathbf{x}_0)$ for an
arbitrarily small $\varepsilon>0$. Let $L_{n1}$ and $L_{n2}$ denote the
distributions of $n^{2/(p+1)} [{\hat{g}}(\mathbf{x}_0)-g(\mathbf{x}_0) ]$ and
$Z_n(\mathbf{0})$, respectively. Then, $d (L_{n1}, L_{n2} )
\rightarrow0$ as $n$ tend to infinity.
\end{theorem}

Computation of the distribution of $Z_n$ solely depends on knowledge of
$\kappa$. Thus one can approximate the distribution of ${\hat{g}}(\mathbf{x}_0)$ by
estimating $\kappa$ and then simulating $Z_n$ with the estimated $\kappa$.
The approximation enables one to correct the downward bias of ${\hat{g}}(\mathbf{x}_0)$
and get an improved estimator of $g(\mathbf{x}_0)$. Estimation of $\kappa$ and
bias-correction for ${\hat{g}}(\mathbf{x}_0)$ will be discussed in Section \ref{sec4}.
\begin{pf*}{Proof of Theorem \protect\ref{theorem2}}
We first give a geometric description of the estimator~${\hat{g}}$.
Consider a hyperplane in $\mathbb{R}^{p}$ defined by
%
%
\begin{equation}\label{sec}
\mathcal{P}(\mathbf{x}_0) = \{\mathbf{x}\in\mathbb{R}_+^{p} \dvtx\mathbf{x}_0^\top(\mathbf{x}-\mathbf{x}_0) = 0\}.
\end{equation}
This hyperplane is perpendicular to the vector $\mathbf{x}_0$ and passes through
$\mathbf{x}_0$. Let $\mathbf{P}_i$ be the point where the ray $\mathbf{R}_i$ meets the hyperplane
$\mathcal{P}^\dag(\mathbf{x}_0) \equiv\mathcal{P}(\mathbf{x}_0) \times\mathbb{R}_+$ in $\mathbb{R}^{p+1}$. It
follows that
%
%
\begin{equation}\label{sec-p}
\mathbf{P}_i = \frac{\|\mathbf{x}_0\|^2}{\mathbf{x}_0^\top\mathbf{X}_i} (\mathbf{X}_i, Y_i ).
\end{equation}
Define $\widehat{\Psi}(\mathbf{x}_0)$ to be the convex-hull of the points $\mathbf{P}
_i$. We claim that
%
%
\begin{equation}\label{prop2}
\widehat{\Psi}(\mathbf{x}_0) = \mathcal{P}^\dag(\mathbf{x}_0) \cap\widehat{\Psi}.
\end{equation}
This means that $\widehat{\Psi}(\mathbf{x}_0)$ is a section of
$\widehat{\Psi}$ obtained by cutting $\widehat{\Psi}$ by the hyperplane
$\mathcal{P}^\dag(\mathbf{x} _0)$. The fact that
$\widehat{\Psi}(\mathbf{x}_0) \subset\mathcal{P}^\dag(\mathbf{x}_0)
\cap \widehat{\Psi}$ follows from convexity of
$\mathcal{P}^\dag(\mathbf{x}_0)$ and~$\widehat{\Psi}$. The reverse
inclusion also holds. To see this, let $(\mathbf{x},y)
\in\mathcal{P}^\dag(\mathbf{x}_0) \cap\widehat{\Psi}$. Since $\widehat
{\Psi}$ is the convex-hull of the rays $\mathbf{R}_i$, it follows that
there exist $\gamma_i^* \ge0$ such that $\mathbf{x}=\sum_{i=1}^n
\gamma_i^* \mathbf{X}_i $ and $y=\sum_{i=1}^n \gamma_i^* Y_i$. Since
$(\mathbf{x},y) \in\mathcal{P}^\dag(\mathbf{x} _0)$, we have
%
%
\begin{equation}\label{prop2-1}
\sum_{i=1}^n \gamma_i^* \mathbf{x}_0^\top\mathbf{X}_i = \|\mathbf{x}_0\|^2.
\end{equation}
Let $\xi_i =
(\mathbf{x}_0^\top\mathbf{X}_i/\|\mathbf{x}_0\|^2)\gamma_i^* \ge0$ for
$1 \le i \le n$. By (\ref{prop2-1}), $\sum_{i=1}^{n+1} \xi_i =1$. By
(\ref {sec-p}), we get $(\mathbf{x},y) = \sum_{i=1}^n \xi_i
\mathbf{P}_i$ which shows $(\mathbf{x},y)
\in\widehat{\Psi}(\mathbf{x}_0)$.

Since $\bigcup_{a \ge0} a \mathcal{P}^\dag(\mathbf{x}_0) =
\mathbb{R}_+^{p+1}$, the CRS property of $\widehat{\Psi}$ and
(\ref{prop2}) thus yield
%
%
\begin{equation}\label{prop3}
\widehat{\Psi} = \bigcup_{a \ge0} a \widehat{\Psi}(\mathbf{x}_0) = \{(a\mathbf{x},
ay)\dvtx
(\mathbf{x}, y) \in\widehat{\Psi}(\mathbf{x}_0), a \ge0\}.
\end{equation}
Recall the definition of ${\hat{g}}$ in (\ref{g-est}). Also, note that,
for $\mathbf{x} \in\mathcal{P}(\mathbf{x}_0)$, we have $(\mathbf{x},y)
\in\widehat{\Psi}$ if and only if $(\mathbf{x},y)
\in\widehat{\Psi}(\mathbf{x}_0)$. This follows from (\ref{prop3}) and
the fact that $a=1$ is the only constant $a \ge0$ such that
$(\mathbf{x},y) \in a \widehat{\Psi}(\mathbf{x}_0)$ if
$\mathbf{x}\in\mathcal{P}(\mathbf{x}_0)$. This gives
%
%
\begin{equation}\label{est-g}
{\hat{g}}(\mathbf{x}) = \sup\{y \dvtx(\mathbf{x},y)
\in\widehat{\Psi}(\mathbf{x}_0)\} \qquad\mbox{if }
\mathbf{x}\in\mathcal{P}(\mathbf{x}_0).
\end{equation}

%
%
\begin{figure}[b]

\includegraphics{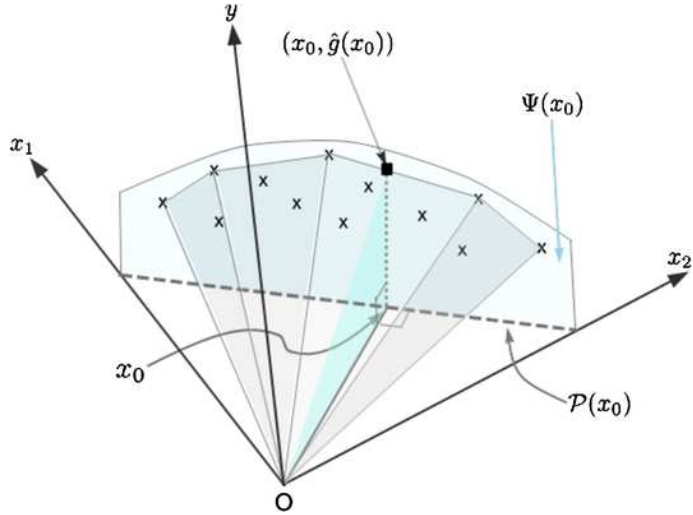}

\caption{An illustration of $\mathcal{P}(\mathbf{x}_0)$, $\mathbf{P}_i$,
$\hat\Psi$ and $\hat{g}$ in the case of $p=2$ and $q=1$. The crosses
are the points $\mathbf{P}_i$, and the gray surface is the roof of the
conical-hull estimator $\hat\Psi$.}\label{fig1}
\end{figure}
%

See Figure \ref{fig1} for an illustration in the case of $p=2$ and
$q=1$.

Let $Q$ be the matrix defined in the paragraph that contains the
definition of the transformation $T_1$ early in this section. Since
$\mathcal{P}(\mathbf{x}_0)= \{\mathbf{x}_0+Q\mathbf{z}_2
\in\mathbb{R}_+^{p} \dvtx\mathbf{z}_2 \in\mathbb{R} ^{p-1}\} $, the
set,
%
%
\begin{eqnarray} \label{sec-1}\qquad
&& \Psi(\mathbf{x}_0) \equiv\{(\mathbf{x}_0+Q\mathbf{z}_2,y) \in A_c
\times\mathbb{R}_+ \dvtx\mathbf{z}_2 \in \mathbb{R}^{p-1}, 0 \le y \le
g(\mathbf{x}_0+Q\mathbf{z}_2)\},
\end{eqnarray}
equals the section of $\Psi$ obtained by cutting $\Psi$ by the
hyperplane $\mathcal{P}^\dag(\mathbf{x}_0)$; that is,
$\Psi(\mathbf{x}_0) = \mathcal{P}^\dag(\mathbf{x}_0) \cap\Psi$. In the
new coordinate system
\[
(\mathbf{z},y')
\equiv\bigl(T_2(\mathbf{x}),y\|\mathbf{x}_0\|^2/(\mathbf{x}_0^\top\mathbf{x})
\bigr),
\]
the set $\Psi(\mathbf{x}_0)$ in (\ref{sec-1}) can be represented by $\{0\}
\times\Psi^*(\mathbf{x}_0)$ where
%
%
\begin{equation}\label{psi-star}
\Psi^*(\mathbf{x}_0) = \{(\mathbf{z}_2,y') \dvtx\mathbf{z}_2
\in\mathbb{R}^{p-1}(\mathbf{x}_0), 0 \le y' \le g^*(\mathbf{z}_2)\}
\end{equation}
and $\mathbb{R}^{p-1}(\mathbf{x}_0)$ denote the set of $\mathbf{z}_2$
such that $\mathbf{x}_0 + Q\mathbf{z}_2 \in A_c$. Also, in that new
coordinate system the points $\mathbf{P}_i$ defined in (\ref{sec-p})
correspond to $(0, \mathbf{P}_i^*)$ where $\mathbf{P}_i^* =
(\mathbf{Z}_{2i}, Y_i')$, $\mathbf{Z}_{2i} =
(\|\mathbf{x}_0\|^2/\mathbf{x}_0^\top\mathbf{X}_i) Q^\top \mathbf{X}_i$
and $Y_i'=(\|\mathbf{x}_0\|^2/\mathbf{x}_0^\top\mathbf{X}_i)Y_i$. Since
convex-hulls are equivariant under linear transformations, this means
that in the new coordinate system, $\widehat{\Psi}(\mathbf{x}_0)$
corresponds to $\{0\} \times \widehat{\Psi}^*(\mathbf{x}_0)$ where
$\widehat{\Psi}^*(\mathbf{x}_0)$ is the convex-hull of the points
$\mathbf{P}_i^*$. Now define
\[
{\hat{g}}^*(\mathbf{z}_2) = {\hat{g}}(\mathbf{x}_0 + Q\mathbf{z}_2)
\]
on $\mathbb{R}^{p-1}(\mathbf{x}_0)$. Since
$(\mathbf{x}_0+Q\mathbf{z}_2, y) \in\widehat{\Psi}(\mathbf{x}_0)$ is
equivalent to $(\mathbf{z}_2,y) \in\widehat{\Psi}^*(\mathbf{x}_0)$, it
follows from (\ref{est-g}) that
%
%
\begin{equation}\label{est-g-star}
{\hat{g}}^*(\mathbf{z}_2) = \sup\{y \dvtx(\mathbf{z}_2,y)
\in\widehat{\Psi}^*(\mathbf{x}_0), \mathbf{z}_2 \in\mathbb{R}^{p-1}\}.
\end{equation}

Let $f$ denote the density of the original random vector $(\mathbf{X},
Y)$ and $f^*$ denote the density of the transformed vector
$(\mathbf{Z}_2, Y')$. The arguments in the preceding paragraph imply
that the distribution of ${\hat{g}}(\mathbf{x}_0)-g(\mathbf{x}_0)$
equals that of ${\hat{g}}^*(\mathbf{0})-g^*(\mathbf{0})$ where
${\hat{g}}^*$ is the convex-hull estimator of $g^*$ constructed from a
random sample of size $n$ generated from the density $f^*$. Let
$\kappa^* = \det(\Lambda)^{-1/2}f^*(\mathbf{0},g^*(\mathbf{0}))$ where
$\Lambda$ is the diagonal matrix with its entries being the eigenvalues
of $-\ddot{g}^*(\mathbf{0})/2$. Define $Z_n^*$ as a version of
${\hat{g}}^*$ constructed from a random sample from the uniform
distribution on $R_n(\kappa^*) \subset\mathbb{R}^p$ where $R_n$ is
defined immediately after (\ref{kappa}). Then one can proceed as in the
proof of Theorem 1 of Jeong and Park (\citeyear{JP06}) to show that the
asymptotic distribution of $n^{2/(p+1)}
({\hat{g}}^*(\mathbf{0})-g^*(\mathbf{0}) )$ is identical to that of
$Z_n^*(\mathbf{0})$ where one uses the transformations $T_3^*$ and
$T_4^*$ defined by
\begin{eqnarray*}
T_3^* \dvtx\mathbf{z}_2 &\mapsto& n^{1/(p+1)}\Lambda^{1/2}P^\top\mathbf{z}_2,\\
T_4^* \dvtx(\mathbf{z}_2,y') &\mapsto& n^{2/(p+1)} \bigl(y'- g^*(\mathbf{0}) -
\dot{g}^*(\mathbf{0})^\top\mathbf{z}_2 \bigr).
\end{eqnarray*}
Recalling the definitions of the transformations $T_3$ and $T_4$ in
(\ref{T3}) and (\ref{T4}), respectively, $T_3^*(\mathbf{z}_2)$ equals $T_3(\mathbf{z})$
without the first component, where $\mathbf{z}^\top=(z_1,\mathbf{z}_2^\top)$, and
$T_4^* (\mathbf{z}_2, y\|\mathbf{x}_0\|/(z_1+\|\mathbf{x}_0\|) ) = T_4(\mathbf{z},y)$. Below, we
prove that $\kappa^*$ equals $\kappa$ defined in (\ref{kappa}) so that
$Z_n^* = Z_n$ in distribution which concludes the proof of the theorem.

Let $T^*$ denote the transformation that maps $(\mathbf{x},y)$ to
\[
(\mathbf{z},y') = \bigl(T_2(\mathbf{x}),
y\|\mathbf{x}_0\|^2/(\mathbf{x}_0^\top\mathbf{x}) \bigr).
\]
Let $c(z_1) = (z_1
+ \|\mathbf{x}_0\|)/\|\mathbf{x}_0\|$. The Jacobian of the inverse transform of $T^*$
equals
\begin{eqnarray*}
J(\mathbf{z}) &\equiv& c(z_1)\det[\|\mathbf{x}_0\|^{-1}(\mathbf{x}_0 +
Q\mathbf{z}_2),
c(z_1) Q ]\\
&=& c(z_1)\,{{\det}^{1/2} }\left[ \matrix{ 1
+(\|\mathbf{z}_2\|/\|\mathbf{x}_0\|)^2&
\bigl(c(z_1)/\|\mathbf{x}_0\|\bigr)\mathbf{z}_2^\top\cr
\bigl(c(z_1)/\|\mathbf{x}_0\|\bigr)\mathbf{z}_2 &
c(z_1)^2 I_{p-1}} \right],
\end{eqnarray*}
where $I_{p-1}$ denotes the identity matrix of dimension $(p-1)$. The
second equality in the above calculation follows from the fact that the
columns of $Q$ are perpendicular to $\mathbf{x}_0$. Thus the joint density of
$T^*(\mathbf{X},Y)$ at the point $(\mathbf{z},y')$ is given by
$J(\mathbf{z})f (c(z_1)(\mathbf{x}_0+Q\mathbf{z}_2), c(z_1)y' )$. The density
$f^*(\mathbf{z}_2,y')$ is simply the marginalization of this joint density with
respect to $z_1$ so that
\[
f^*(\mathbf{z}_2,y') = \int_{-\|\mathbf{x}_0\|}^\infty
J(\mathbf{z})f \bigl(c(z_1)(\mathbf{x}_0+Q\mathbf{z}_2), c(z_1)y' \bigr)\,dz_1.
\]
Now, since $J(z_1,\mathbf{0})=c(z_1)^p$, we obtain
\begin{eqnarray*}
f^*(\mathbf{0},g^*(\mathbf{0})) &=& \int_{-\|\mathbf{x}_0\|}^\infty c(z_1)^p
f (c(z_1)\mathbf{x}_0, c(z_1)g^*(\mathbf{0}) ) \,dz_1\\
&=& \theta,
\end{eqnarray*}
where $\theta$ is defined in (\ref{theta}).
\end{pf*}

\begin{figure}

\includegraphics{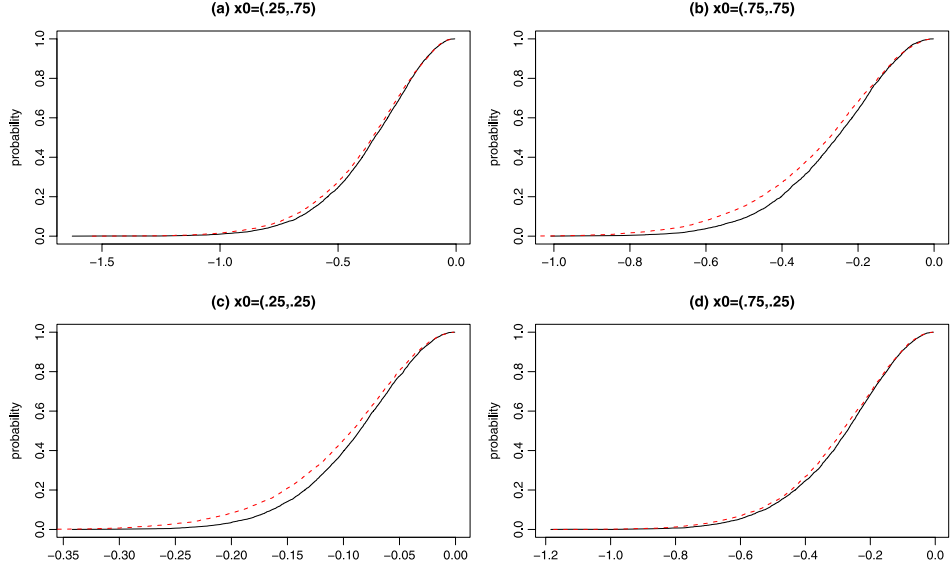}

\caption{Solid curves are the empirical distribution
functions of $Z_n(0)$, and the dotted curves are those of
$n^{2/(p+1)} \{{\hat{g}}(\mathbf{x}_0)-g(\mathbf{x}_0) \}$ in the case where
$n=100$ and $\lambda=3$.}\vspace*{-5pt}
\label{fig2}
\end{figure}

\begin{figure}[t]\vspace*{-5pt}

\includegraphics{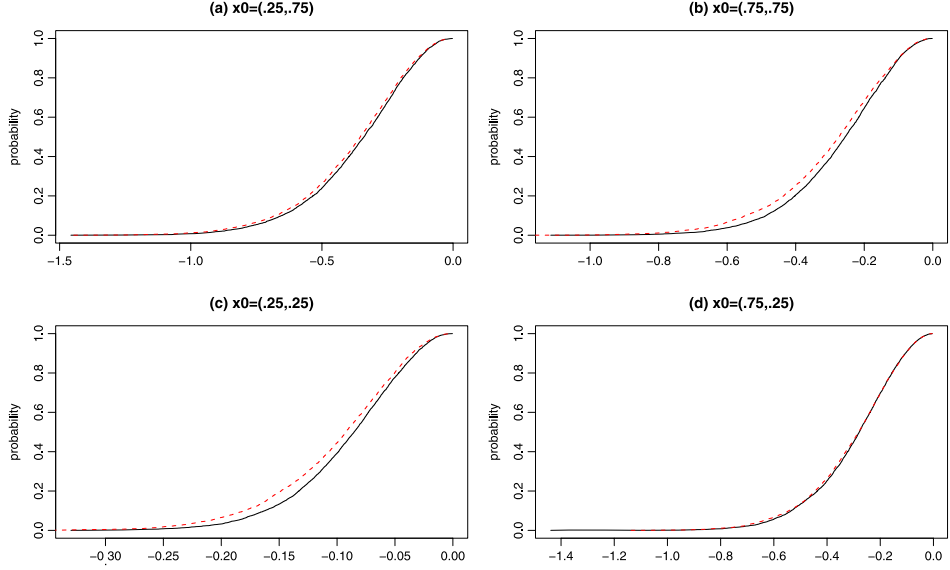}

\caption{Solid curves are the empirical distribution
functions of $Z_n(0)$, and the dotted curves are those of
$n^{2/(p+1)} \{{\hat{g}}(\mathbf{x}_0)-g(\mathbf{x}_0) \}$ in the case where
$n=400$ and $\lambda=3$.}\label{fig3}
\end{figure}

To see how well the distribution of $n^{2/(p+1)} \{{\hat{g}}(\mathbf{x}_0)-g(\mathbf{x}
_0) \}$ is approximated by
that of $Z_{n}(\mathbf{0})$, we took a Cobb--Douglas CRS production
function $g(\mathbf{x})=x_1^{0.4}\times x_2^{0.6}$ ($p=2$). We generated 5000
random samples of size $n=100$ and $400$ from $f(x_1,x_2,y)=\lambda
x_1^{-0.4\lambda}x_2^{-0.6\lambda} y^{\lambda-1}$ supported on $\Psi=\{
(x_1,x_2,y) \dvtx0\le x_1,x_2\le1, 0\le y\le g(x_1,x_2)\}$. This
yielded i.i.d. copies of $(X_1,X_2,Y)$ with $X_1\sim
\operatorname{Uniform}[0,1]$, $X_2\sim\operatorname{Uniform}[0,1]$ and
$Y=\break g(X_1,X_2)e^{-V/\lambda}$ where $V\sim\operatorname{Exp}(1)$.
Figures \ref
{fig2} and \ref{fig3} depict the empirical distributions of
$n^{2/(p+1)} \{{\hat{g}}(\mathbf{x}_0)-g(\mathbf{x}_0) \}$ and $Z_{n}(\mathbf{0})$
based on these samples in the case where $\lambda=3$. The figures
suggest that the approximation is fairly good for moderate sample sizes
and get better as the sample size increases.

%

%

Theorem \ref{theorem2} excludes the case where $g$ is linear; that is, $g(\mathbf{x})=\mathbf{c}
^\top\mathbf{x}$ for some vector $\mathbf{c}$. The latter case needs a different
treatment. In the following theorem, we give the limit distribution in
this case. To state the theorem, let $(\mathbf{V}_{2i}^L, W_i^L)$ be a random
sample from the uniform distribution on the $p$-dimensional rectangle,
%
%
\begin{eqnarray} \label{rect}
R_n^L(\theta) &=& \bigl[ -\tfrac{1}{2}\theta^{-1/(p+1)} n^{1/(p+1)},
\tfrac{1}{2}\theta^{-1/(p+1)} n^{1/(p+1)}
\bigr]^{p-1}\nonumber\\[-8pt]\\[-8pt]
&&{} \times\bigl[ -\theta^{-2/(p+1)}n^{2/(p+q)}, 0\bigr],\nonumber
\end{eqnarray}
where $\theta$ is defined in (\ref{theta}). The volume of this set in
$\mathbb{R}^p$ equals $n \theta^{-1}$. Let $Z_n^L(\cdot)$ be a version of
$Z_n(\cdot)$ constructed from $(\mathbf{V}_{2i}^L, W_i^L)$ replacing $(\mathbf{V}_{2i},
W_i)$.
\begin{theorem}\label{theorem3}
Assume \textup{(A1)} and \textup{(A2)}. Assume further that $\Psi=
\{(\mathbf{x},y) \in\mathbb{R}_+^{p+1} \dvtx0 \le y
\le\mathbf{c}^\top\mathbf{x}\}$ for some constant vector
$\mathbf{c}\neq\mathbf{0}$ and that the density $f$ of $(\mathbf{X},
Y)$ is uniformly continuous on $H_\varepsilon(\mathbf{x}_0)$ for an
arbitrarily small $\varepsilon>0$. Let $L_{n1}$ and $L_{n2}'$ denote
the distributions of $n^{2/(p+1)}
[{\hat{g}}(\mathbf{x}_0)-\mathbf{c}^\top\mathbf{x}_0 ]$ and
$Z_n^L(\mathbf{0} )$, respectively. Then $d (L_{n1}, L_{n2}' )
\rightarrow0$ as $n$ tends to infinity.
\end{theorem}
\begin{pf} In this case we consider the following
transformation:
%
%
\begin{equation}\label{trans-l}
T^L \dvtx(\mathbf{x}, y) \mapsto\bigl(T_3^L \circ T_2 (\mathbf{x}),
T_4^L (T_2 (\mathbf{x}),y ) \bigr),
\end{equation}
where $T_3^L \dvtx\mathbf{z}\mapsto(z_1, n^{1/(p+1)}\mathbf{z}_2^\top)^\top$
and
\[
T_4^L \dvtx(\mathbf{z},y) \mapsto n^{2/(p+2)}
\biggl(\frac{\|\mathbf{x}_0\|}{z_1+\|\mathbf{x}_0\|}y -
\mathbf{c}^\top\mathbf{x}_0 - \mathbf{c}^\top Q \mathbf{z}_2 \biggr).
\]
Let $(\mathbf{V}^L,W^L) = T^L(\mathbf{X},Y)$. Then it can be shown as in the proof of
Theorem \ref{theorem2} that the density of $(\mathbf{V}_2^L,W^L)$ is given by $n^{-1}\theta
\{1+o(1)\}$ uniformly for $\mathbf{v}_2^L$ and $w^L$ in any compact sets of
respective dimension. The rest of the proof is the same as that for
Theorem \ref{theorem2}.
\end{pf}

In the special case where $p=1$, we can derive the limit distribution
explicitly. In this case, the boundary function $g$ is linear and takes
the form $g(x)=cx$ for some constant $c>0$. The transformation $T^L$ in
(\ref{trans-l}) reduces to
\[
T^L(x,y) = \biggl(x-x_0, n \biggl(\frac{y}{x}x_0 -cx_0 \biggr) \biggr).
\]
The marginal density of $W^L$, where $(V^L,W^L)=T^L(X,Y)$, is
approximated by the constant $n^{-1}\theta$ uniformly for $w^L$ in any
compact subset of $\mathbb{R}_-$ where $\theta$ in this case equals
$x_0 \int_0^\infty u f(ux_0, ucx_0) \, du$. According to Theorem
\ref{theorem3}, the limit distribution of $n({\hat{g}}(x_0)-g(x_0))$
equals the limit distribution of $Z_n^L$ which is nothing else than
$\max_{i=1}^n W_i^L$ in this simplest case where $W_i^L$ are a random
sample from the uniform distribution on $[-n\theta^{-1}, 0]$. Since
$-\max_{i=1}^n W_i^L$ has the exponential distribution with mean
$\theta^{-1}$ in the limit, we have
\[
P \bigl[n \bigl(g(x_0)-{\hat{g}}(x_0) \bigr) \le w \bigr] \rightarrow
1-\exp(-\theta w)
\]
for all $w \ge0$.

\subsection{The case where $q>1$}\label{sec32}

In this section we extend the results in the previous section to the
case where $q > 1$ and $\Psi$ is a conical-hull of a convex set
$\mathcal{A}$ in $\mathbb{R} _+^{p+q}$. For this we make a canonical
transformation on $\mathbf{y}$-space so that the problem for $q>1$ is
reduced to the case where $q=1$. Again we fix the point
$(\mathbf{x}_0,\mathbf{y}_0)$ where we want to estimate the function
$\lambda$.

Let $\Gamma$ be a $q \times(q-1)$ matrix whose columns form a basis for
$\mathbf{y}_0^\perp$. Consider a transformation $\mathcal{T}$ that maps
$\mathbf{y}\in\mathbb{R}_+^q$ to $(\mathbf{u},\omega)
\in\mathbb{R}^{q-1} \times\mathbb{R}_+$ where
%
\begin{equation}\label{roty0}
\mathbf{u}= \Gamma^\top\mathbf{y},\qquad
\omega=\frac{\mathbf{y}_0^\top\mathbf{y}}{\|\mathbf{y}_0\|}.
\end{equation}
Then, in the new coordinate system $(\mathbf{x},\mathbf{u},\omega)$,
the set $\Psi$ can be represented as
%
%
\begin{equation}\label{psi-t}
\Psi_\mathcal{T}= \biggl\{(\mathbf{x},\mathbf{u},\omega) \in\mathbb{R}_+^p
\times\mathbb{R}^{q-1} \times \mathbb{R}_+ \dvtx\biggl(\mathbf{x},
\Gamma\mathbf{u}+ \omega \frac{\mathbf{y}_0}{\|\mathbf{y}_0\|} \biggr)
\in\Psi\biggr\}.
\end{equation}
Define a $(p+q-1)$-dimensional function
\[
g_\mathcal{T}(\mathbf{x},\mathbf{u}) \equiv
g_\mathcal{T}(\mathbf{x},\mathbf{u};\mathbf{y}_0) = \sup\biggl\{a>0
\dvtx\biggl(\mathbf{x}, \Gamma\mathbf{u}+ a
\frac{\mathbf{y}_0}{\|\mathbf{y}_0\|} \biggr) \in\Psi\biggr\}.
\]
This is a boundary function in the transformed space such that all
points $(\mathbf{x},\mathbf{u},\omega)$ in $\Psi_\mathcal{T}$ lie below
the surface represented by the equation
$\omega=g(\mathbf{x},\mathbf{u})$.

Convexity of the function $g_\mathcal{T}$ follows from the fact that,
due to convexity of $\Psi$,
\[
a_0 \in\biggl\{a>0 \dvtx\biggl(\mathbf{x}, \Gamma\mathbf{u}+ a
\frac{\mathbf{y}_0}{\|\mathbf{y}_0\|} \biggr) \in\Psi\biggr\}
\]
and
\[
a_0' \in\biggl\{a'>0 \dvtx\biggl(\mathbf{x}', \Gamma\mathbf{u}' + a' \frac{\mathbf{y}_0}{\|\mathbf{y}
_0\|} \biggr) \in\Psi\biggr\},
\]
together, imply
\begin{eqnarray*}
&& \alpha a_0 + (1-\alpha)a_0'
\\
&&\qquad\in\biggl\{a>0
\dvtx\biggl(\alpha\mathbf{x}+(1-\alpha)\mathbf{x}',
\Gamma\bigl(\alpha\mathbf{u}+(1-\alpha)\mathbf{u}'\bigr)+ a
\frac{\mathbf{y}_0}{\|\mathbf{y}_0\|} \biggr) \in\Psi\biggr\}.
\end{eqnarray*}
Also, it has the CRS property (\ref{crs-1}) since $\Psi$ satisfies
(\ref {crs}). Furthermore, since $(\mathbf{x},\mathbf{y}) \in\Psi$ if
and only if $(\mathbf{x}, \mathcal{T}(\mathbf{y}))
\in\Psi_\mathcal{T}$, and $\mathcal{T}(\alpha\mathbf{y}_0) =
(\mathbf{0}^\top,\alpha\|\mathbf{y}_0\|)^\top$ for all $\alpha>0$, we
obtain
%
%
\begin{eqnarray}\label{eff-equiv}
g_\mathcal{T}(\mathbf{x}_0,\mathbf{0}) &=& \sup\biggl\{a>0 \dvtx\biggl(\mathbf{x}_0, a
\frac{\mathbf{y}_0}{\|\mathbf{y}_0\|} \biggr) \in\Psi\biggr\} \nonumber\\
&=& \sup\{a>0 \dvtx(\mathbf{x}_0, (\mathbf{0},a)) \in
\Psi_\mathcal{T}\} \nonumber\\
&=& \|\mathbf{y}_0\| \sup\{\lambda>0 \dvtx(\mathbf{x}_0, (\mathbf{0},\lambda
\|\mathbf{y}_0\|) )
\in\Psi_\mathcal{T}\} \\
&=& \|\mathbf{y}_0\| \sup\{\lambda>0 \dvtx(\mathbf{x}_0, \mathcal{T}(\lambda
\mathbf{y}_0) ) \in\Psi_\mathcal{T}\} \nonumber\\
\label{gT}
&=& \|\mathbf{y}_0\| \lambda(\mathbf{x}_0,\mathbf{y}_0).\nonumber
\end{eqnarray}
Here and below, $\mathbf{0}$ denotes the $(q-1)$-dimensional zero
vector. Thus the problem of estimating
$\lambda(\mathbf{x}_0,\mathbf{y}_0)$ using
$(\mathbf{X}_i,\mathbf{Y}_i)$ is reduced to that of estimating
$g_\mathcal{T}(\mathbf{x}_0,\mathbf{0})$ in the transformed space using
$(\mathbf{X}_i, \mathcal{T}(\mathbf{Y}_i))$.

We note that in the proof of Theorem \ref{theorem2} we use only
convexity and the CRS property of $g$. Thus the theory we developed in
the previous section is applicable to $g_\mathcal{T}$. Let
$(\mathbf{U}_i, \Omega_i) = \mathcal{T}(\mathbf{Y}_i)$ where
$\mathbf{U}_i$ is the vector of the first $(q-1)$ elements of
$\mathcal{T}(\mathbf{Y}_i)$, and $\Omega_i$ is the scalar-valued random
variable. The joint density of $(\mathbf{X}_i,\mathbf{U}_i, \Omega_i)$
at the point $(\mathbf{x},\mathbf{u},\omega)$ is given by
%
%
\begin{equation}\label{trans-den}
f_\mathcal{T}(\mathbf{x},\mathbf{u},\omega) =
{\det}^{1/2}(\Gamma^\top\Gamma)f \biggl(\mathbf{x}, \Gamma\mathbf{u}+
\omega\frac{\mathbf{y}_0}{\|\mathbf{y}_0\|} \biggr).
\end{equation}
The constant $\theta$ defined in (\ref{theta}) that corresponds to the
density $f_\mathcal{T}$ equals
\begin{eqnarray*}
\theta_\mathcal{T}&=& \|(\mathbf{x}_0,\mathbf{0})\| \int_0^\infty
u^{p+q-1} f_\mathcal{T}
(u\mathbf{x}_0,\mathbf{0}, u g_\mathcal{T}(\mathbf{x}_0,\mathbf{0}) ) \,du\\
&=&{\det}^{1/2}(\Gamma^\top\Gamma)\|\mathbf{x}_0\| \int_0^\infty
u^{p+q-1}f \biggl(u\mathbf{x}_0,u g_\mathcal{T}
(\mathbf{x}_0,\mathbf{0})\frac{\mathbf{y}_0}{\|\mathbf{y}_0\|} \biggr) \,du\\
&=& {\det}^{1/2}(\Gamma^\top\Gamma)\|\mathbf{x}_0\| \int_0^\infty
u^{p+q-1} f (u\mathbf{x}_0,u
\lambda(\mathbf{x}_0,\mathbf{y}_0)\mathbf{y}_0 ) \,du,
\end{eqnarray*}
where the last identity follows from (\ref{eff-equiv}). The determinant
that corresponds to $\det(\Lambda)$ in the definition of $\kappa$ in
(\ref{kappa}) is $\det(-Q_\mathcal{T}^\top
\ddot{g}_\mathcal{T}(\mathbf{x}_0,\mathbf{0})Q_\mathcal{T}/2)$ where
$Q_\mathcal{T}$ is a $(p+q-1) \times (p+q-2)$ matrix whose columns form
an orthonormal basis for $(\mathbf{x}_0,\mathbf{0})^\perp$. Thus we
modify the definition of $\kappa$ as
\[
\kappa_\mathcal{T}= \theta_\mathcal{T}\det\bigl(-Q_\mathcal{T}^\top
\ddot{g}_\mathcal{T}(\mathbf{x}_0,\mathbf{0})Q_\mathcal{T}/2\bigr)^{-1/2}.
\]
Recall that the construction of $Z_n$ defined in (\ref{limit}) depends
only on $\kappa$ and $p$. Define $Z_{n,\mathcal{T}}$ as a version of
$Z_n$ with $\kappa_\mathcal{T}$ and $(p+q-1)$ replacing $\kappa$ and
$p$, respectively. Also, define a $(p+q-2)$-dimensional function
$g_\mathcal{T}^*(\mathbf{z}_2) = g_\mathcal{T}((\mathbf{x}
_0,\mathbf{0})+Q_\mathcal{T}\mathbf{z}_2)$, and
$H_{\varepsilon,\mathcal{T}}(\mathbf{x}_0,\mathbf{0})$ as
$H_{\varepsilon}(\mathbf{x}_0)$ at (\ref{neigh}) with $(p+q-1)$,
$g_\mathcal{T}$, $(\mathbf{x} _0, \mathbf{0})$ and $Q_\mathcal{T}$
replacing $p$, $g$, $\mathbf{x}_0$ and $Q$, respectively. Then we have
the following theorem for the limit distribution of
${\hat{\lambda}}(\mathbf{x}_0,\mathbf{y}_0)$ for arbitrary dimensions
$p, q \ge1$.
\begin{theorem}\label{theorem4}
Assume \textup{(A1)} and \textup{(A2)}. In addition, assume that
$-\ddot{g}_\mathcal{T} ^*$ is positive definite and continuous at
$\mathbf{0}$, and that the density $f_\mathcal{T}$ given at
(\ref{trans-den}) is uniformly continuous on
$H_{\varepsilon,\mathcal{T}}(\mathbf{x}_0,\mathbf{0})$ for an
arbitrarily small $\varepsilon>0$. Let $L_{n1}$ and $L_{n2}$ denote the
distributions of $n^{2/(p+q)}
[{\hat{\lambda}}(\mathbf{x}_0,\mathbf{y}_0)-\lambda(\mathbf{x}_0,\mathbf{y}_0)
]$ and $Z_{n,\mathcal{T}}(\mathbf{0}_{p+q-2})/\|\mathbf{y}_0\|$,
respectively. Then, $d (L_{n1}, L_{n2} ) \rightarrow0$ as $n$ tends to
infinity.
\end{theorem}

Theorem \ref{theorem4} excludes the case where $\Psi=
\{(\mathbf{x},\mathbf{y}) \in\mathbb{R}_+^{p+q}\dvtx
\mathbf{c}_1^\top\mathbf{x}- \mathbf{c}_2^\top\mathbf{y}\ge0\}$ for
some constant vectors $\mathbf{c}_1, \mathbf{c}_2 > \mathbf{0}$. Below
we treat this case. When $q=1$, this corresponds to the case where the
boundary function $g$ is linear in $\mathbf{x}$.

Define
\[
\mathbf{c}_\mathcal{T}=
\frac{\|\mathbf{y}_0\|}{\mathbf{c}_2^\top\mathbf{y}_0}
\pmatrix{\mathbf{c}_1 \cr \Gamma^\top(-\mathbf{c}_2)}.
\]
Then $\Psi_\mathcal{T}$ defined in (\ref{psi-t}) takes the form
\[
\Psi_\mathcal{T}= \left\{(\mathbf{x},\mathbf{u},w) \dvtx0 \le w
\le\mathbf{c}_\mathcal{T}^\top\pmatrix{\mathbf{x} \cr \mathbf{u}} \right\},
\]
and it holds that
\[
\mathbf{c}_\mathcal{T}^\top\pmatrix{\mathbf{x}\cr\mathbf{0}} =
\|\mathbf{y}_0\| \lambda(\mathbf{x}_0,\mathbf{y}_0).
\]
Thus we can apply the arguments leading to Theorem \ref{theorem3} with
$p$, $\mathbf{c}$, $\mathbf{x}_0$ and $Q$ being replaced by $(p+q-1)$,
$\mathbf{c}_\mathcal{T}$, $(\mathbf{x}_0, \mathbf{0})$ and
$Q_\mathcal{T}$, respectively.

Let $R_{n,\mathcal{T}}^L(\theta_cT)$ be the rectangle defined in
(\ref{rect}) with $\theta$ and $p$ being replaced by
$\theta_\mathcal{T}$ and $(p+q-1)$. Define $Z_{n,\mathcal{T}}^L$ as
$Z_n^L$ using a random sample from the uniform distribution of the
$(p+q-1)$-dimensional rectangle $R_{n,\mathcal{T}}^L(\theta
_\mathcal{T})$. By applying the proof of Theorem \ref{theorem3} to
$\mathbf{c}_\mathcal{T}$ replacing $\mathbf{c} $, we get the following
theorem.
\begin{theorem}\label{theorem5}
Assume \textup{(A1)} and \textup{(A2)}. Assume further that $\Psi=
\{(\mathbf{x},\mathbf{y}) \in\mathbb{R}_+^{p+q}\dvtx
\mathbf{c}_1^\top\mathbf{x}- \mathbf{c}_2^\top\mathbf{y}\ge0\}$ for
some constant vectors $\mathbf{c}_1, \mathbf{c}_2 > \mathbf{0}$ and
that the density $f_\mathcal{T}$ given at (\ref {trans-den}) is
uniformly continuous on
$H_{\varepsilon,\mathcal{T}}(\mathbf{x}_0,\mathbf{0})$ for an
arbitrarily small $\varepsilon>0$. Let $L_{n1}$ and $L_{n2}'$ denote
the distributions of $n^{2/(p+q)}
[{\hat{\lambda}}(\mathbf{x}_0,\mathbf{y}_0)-\lambda(\mathbf{x}_0,\mathbf{y}
_0) ]$ and $Z_{n,\mathcal{T}}^L(\mathbf{0}_{p+q-2})/\|\mathbf{y}_0\|$,
respectively. Then $d (L_{n1}, L_{n2}' ) \rightarrow0$ as $n$ tends to
infinity.
\end{theorem}

\section{Estimation of $\kappa$ and $\kappa_\mathcal{T}$}\label{sec4}

We discuss how to estimate $\kappa$ as defined in ({\ref{kappa}) for
the case where $q=1$. It is straightforward to extend the methods to
the case where $q>1$ via the canonical transformation that we introduced
in Section \ref{sec32}.

Consider the set $H_\varepsilon(\mathbf{x}_0) \subset\mathbb{R}_+^{p+1}$ defined in
(\ref{neigh}). The projection of this set on the $\mathbf{x}$-space is a conical
hull around the vector $\mathbf{x}_0$, and for each direction of the ray $\mathbf{x}
_0 +
Q\mathbf{z}_2$, determined by $\mathbf{z}_2$, its section on that direction is also a
conical hull of single dimension under the boundary $g$. For each fixed
$u \ge0$, let
\begin{eqnarray*}
&& H_\varepsilon(u;\mathbf{x}_0) = \bigl\{ \bigl(u(\mathbf{x}_0
+ Q\mathbf{z}_2),y \bigr) \dvtx\|\mathbf{z}_2\|\le\varepsilon,\\
&&\hspace*{62.3pt}
g \bigl(u(\mathbf{x}_0 + Q\mathbf{z}_2) \bigr)-u\varepsilon\le y \le
g \bigl(u(\mathbf{x}_0 + Q\mathbf{z}_2) \bigr) \bigr\}.
\end{eqnarray*}
This is a section of $H_\varepsilon(\mathbf{x}_0)$ obtained by cutting
$H_\varepsilon(\mathbf{x}_0)$ perpendicular to $\mathbf{x}_0$ at the
distance $u \|\mathbf{x} _0\|$ from the origin. Its volume in the
cutting hyperplane $u \mathcal{P}^\dag(\mathbf{x}_0)$, where
$\mathcal{P}^\dag(\mathbf{x}_0)$ is defined between (\ref{sec}) and
(\ref{sec-p}), equals
\[
v_\varepsilon(u) = c_{p-1} u^p \varepsilon^p,
\]
where $c_r$ denote the volume of the $r$-dimensional unit ball, that
is, $c_r=\frac{\pi^{r/2}}{\Gamma(r/2+1)}$ with $\Gamma(z)=\int_0^\infty
t^{z-1}e^{-t} dt$. Thus, as $\varepsilon\rightarrow0$ we have
\begin{eqnarray*}
P [(\mathbf{X},Y) \in H_\varepsilon(\mathbf{x}_0) ] &=& \int_0^\infty
\int_{(\mathbf{x},y) \in H_\varepsilon(u;\mathbf{x}_0)} f(\mathbf{x},y) \,d\mathbf{x} \,dy \,du\\
&=& \int_0^\infty f(u\mathbf{x}_0, u g(\mathbf{x}_0)) v_\varepsilon(u) \,du\, \{1+o(1)\}
\\
&=& c_{p-1} \varepsilon^p \int_0^\infty u^p f(u\mathbf{x}_0, u
g(\mathbf{x}_0)) \,du\, \{1+o(1)\}.
\end{eqnarray*}
This consideration motivates the following estimator of $\theta$:
%
%
\begin{equation}\label{est-th}
{\hat{\theta}}= \|\mathbf{x}_0\| c_{p-1}^{-1} n^{-1}\varepsilon^{-p}
\sum_{i=1}^n I \bigl((\mathbf{X}_i,Y_i)
\in\widehat{H}_\varepsilon(\mathbf{x}_0) \bigr),
\end{equation}
where $\hat{H}_\varepsilon(\mathbf{x}_0)$ is the sample version of
$H_\varepsilon(\mathbf{x} _0)$ with $g$ replaced by $\hat{g}$ in its
definition. Note that, for implementing ${\hat{\theta}}$, it is
convenient to use the fact,
\[
(\mathbf{X}_i,Y_i) \in\hat{H}_\varepsilon(\mathbf{x}_0)
\quad\Leftrightarrow\quad\|\mathbf{Z}_{2i}\|\le\varepsilon,\qquad \hat{g}^*(\mathbf{Z}
_{2i})-\varepsilon\le Y_i^\prime\le\hat{g}^*(\mathbf{Z}_{2i}).
\]
It is straightforward to see that ${\hat{\theta}}$ is a consistent
estimator of $\theta$ under the conditions of Theorem \ref{theorem2}.

For estimating $\det(\Lambda)$, one can apply local polynomial fitting
to $ \{ (\mathbf{Z}_{2i},\break\hat{g}^*(\mathbf{Z}_{2i}) ) \}$. For a small
$\delta>0$, perform a second-order polynomial regression on the set of
the points
\[
\{(\mathbf{Z}_{2i},\hat{g}^*(\mathbf{Z}_{2i}))
\dvtx\|\mathbf{Z}_{2i}\|\le\delta, i=1,2,\ldots,n
\}\cup\{(\mathbf{0},\hat{g}^*(\mathbf{0}) \},
\]
to get
%
\begin{equation}\label{est-lamb}
\breve{g}^*(\mathbf{z})=\breve{g}_0+\breve{\mathbf{g}}_1^\prime\mathbf{z}+{\mathbf{z}}^\prime\breve
{\mathbf{g}}_2\mathbf{z}.
\end{equation}
Use $\det({\breve{\mathbf{g}}_2})$ as an estimator of
$\det({\Lambda})$. An estimator of $\kappa$ is then defined by
$\hat{\kappa} = {\hat{\theta}}\det ({\breve{\mathbf{g}}_2})^{-1/2}$.

Using the estimator of $\kappa$ one can obtain a bias-corrected
estimator of the function $g^*$. For this, one generates $Z_n$
repeatedly as described at (\ref{limit}) using the estimated $\kappa$.
Call them $Z_{n,1},Z_{n,2},\ldots,Z_{n,B}$. A bias-corrected estimator
is then defined by
\[
\hat{g}^*(\mathbf{0})-n^{-2/(p+1)}\bar{Z}_{n,\cdot}(\mathbf{0}),
\]
where $\bar{Z}_{n,\cdot}(\mathbf{0})=B^{-1}\sum_{b=1}^B Z_{n,b}(\mathbf{0})$.
Also, a $100\times(1-\alpha)\%$ confidence interval is given by
\[
\bigl[ \hat{g}^*(\mathbf{0})-n^{-2/(p+1)}Z_{n,(B(1-\alpha/2))}(\mathbf{0}), \hat
{g}^*(\mathbf{0})-n^{-2/(p+1)}Z_{n,(B\alpha/2)}(\mathbf{0})\bigr],
\]
where $Z_{n,(j)}(\mathbf{0})$ are the ordered values $Z_{n,j}(\mathbf{0})$ such
that $Z_{n,(1)}(\mathbf{0}) > Z_{n,(2)}(\mathbf{0}) > \cdots> Z_{n,(B)}(\mathbf{0})$.

\section{Numerical study}\label{sec5}

In this section we investigate, by a Monte Carlo experiment, the
behavior of the sampling distribution of the DEA--CRS estimator in
finite samples. To be more specific we will compare if the
bias-corrected estimator suggested above has better properties than the
original DEA--CRS estimator in terms of median squared error.

For our Monte Carlo scenario, we adapted the scenario proposed in
Kneip, Simar and Wilson (\citeyear{KSW08}) to our setup. The efficient
frontier is
defined with a CRS generalized Cobb--Douglas production function,
\begin{eqnarray*}
Y_{1e}&=&X_{1}^{0.4}X_{2}^{0.6} \cos\omega,\\
Y_{2e}&=&X_{1}^{0.5}X_{2}^{0.5} \sin\omega,
\end{eqnarray*}
where the random rays are generated through $\omega\sim
\operatorname{Uniform}(\frac{1}{9}\frac{\pi}{2},\frac{8}{9}\frac{\pi
}{2})$ and the
values of the inputs $\mathbf{X}$ by $(X_{1},X_{2})\sim\operatorname
{Uniform}[10,20]^2$.
Then inefficient firms are generated below the efficient frontier by
\[
(Y_1,Y_2)=(Y_{1e},Y_{2e}) e^{-V/3}\qquad \mbox{where } V\sim\operatorname{Exp}(1).
\]

So we are in a situation with $p=q=2$, and we will analyze the
estimation of the efficiency score of the fixed point $\mathbf{x}_0=(15,15)$,
$\mathbf{y}_0=(10,10)$. It is easy to see that the true value of the parameter
to estimate is $\lambda_0=\lambda(\mathbf{x}_0, \mathbf{y}_0)=1.0607$. We analyze the
cases $n=100$ and $n=400$.

We performed 500 Monte Carlo simulations and computed the squared
errors of the original DEA--CRS estimator and of the bias-corrected
estimator. Table \ref{sim} summarizes the results. It gives the ratios
of the median of the squared error of the two estimators,
\[
\mathrm{R}_{\varepsilon,\delta}= \frac{\mathrm{med}\{(\tilde\lambda
_{0,j}-\lambda
_0)^2, j=1,2,\ldots,500\}}{\mathrm{med}\{({\hat{\lambda}}_{0,j}-\lambda
_0)^2, j=1,2,\ldots,500\}},
\]
where ${\hat{\lambda}}_{0,j}$ and $\tilde\lambda_{0,j}$ denote the original
DEA--CRS estimate and the bias-corrected estimate computed in the $j$th
Monte Carlo replication, respectively. Note that the bias-corrected
estimator relies on the values of the smoothing parameters $(\varepsilon
,\delta)$ which appear in the definitions (\ref{est-th}) and (\ref
{est-lamb}), respectively.

%
\begin{table}
\caption{Ratio $\mathrm{R}_{\varepsilon,\delta}$ of the median of
the squared errors of the bias-corrected estimator over the median of
the squared errors of the original DEA--CRS estimator}\label{sim}
\begin{tabular*}{\tablewidth}{@{\extracolsep{\fill}}lclc@{}}
\hline
\multicolumn{2}{@{}c}{\textbf{(}$\bolds{n=100}$\textbf{)}} & \multicolumn{2}{c@{}}{\textbf{(}$\bolds{n=400}$\textbf{)}}
\\[-4pt]
\multicolumn{2}{@{}c}{\hrulefill} & \multicolumn{2}{c@{}}{\hrulefill} \\
& \textbf{Ratio of median} & & \textbf{Ratio of median} \\
$\bolds{\varepsilon=\delta}$ & \textbf{of squared errors} &
\multicolumn{1}{c}{$\bolds{\varepsilon=\delta}$}
& \textbf{of squared errors}\\
\hline
3.50 & 0.7123 & 3.25 & 0.6500 \\
3.75 & 0.6863 & 3.50 & 0.6402 \\
4.00 & 0.7264 & 3.75 & 0.6965 \\
4.25 & 0.8081 & 4.00 & 0.7026  \\
4.50 & 0.8213 & 4.25 & 0.7734 \\
\hline
\end{tabular*}
\end{table}

It is observed from the table that the bias-correction works very well
for a wide range of the smoothing parameters, even though the smoothing
parameters were taken to be equal in the simulation study for saving
computational costs. We see also that the performance of the
bias-corrected estimator gets better when compared to the original
DEA--CRS as the sample size increases.

\section{Discussion}\label{sec6}

In this paper we developed the theoretical properties of the DEA
estimator defined in (\ref{dea-est}) in the case where the support $\Psi
$ of the data $(\mathbf{X}_i, \mathbf{Y}_i)$ satisfies the CRS condition (\ref{crs}).
The assumption of CRS may be tested. In fact, whether the underlying
technology exhibits CRS or VRS is a crucial question in studying
productive efficiency. The question has important economic
implications. If the technology does not exhibit CRS, then some
production units may be found to be either too large or too small.
Using the estimator at (\ref{dea-est}) in the case where the true
technology displays nonconstant returns to scale results in
statistically inconsistent estimates of efficiency and seriously
distorts measures of efficiency.

One way to test CRS against VRS is to use the test statistic defined as
\[
\rho_n = \frac{1}{n}\sum_{i=1}^n
\biggl(\frac{{\hat{\lambda}}(\mathbf{X}_i,\mathbf{Y}_i)}{{\hat{\lambda}}
_{\mathrm{VRS}}(\mathbf{X}_i,\mathbf{Y}_i)} -1 \biggr),
\]
where ${\hat{\lambda}}_{\mathrm{VRS}}$ is a version of ${\hat{\lambda}}$ for the case of VRS
defined as in (\ref{dea-est}) but with $\widehat{\Psi}$ replaced by the
convex-hull of $\{(\mathbf{X}_i, \mathbf{Y}_i)\}_{i=1}^n$. By construction,
\[
{\hat{\lambda}}(\mathbf{X}_i,\mathbf{Y}_i) \ge{\hat{\lambda}}_{\mathrm{VRS}}(\mathbf{X}_i,\mathbf{Y}_i) > 0
\]
so that $\rho_n \ge0$. A larger value of $\rho_n$ gives a stronger
evidence against the null hypothesis of CRS in favor of the alternative
hypothesis of VRS. The test statistic was considered by Simar and
Wilson (\citeyear{SW02}). One may compute $p$-values or critical values
using a
bootstrap method. For example, a subsampling scheme with the subsample
size determined by the procedure described in Politis, Romano and Wolf
(\citeyear{PRW01}) might work for this problem. For testing CRS against
nonconstant
returns-to-scale, which is broader than VRS, one may use the estimators
analyzed by Hall, Park and Stern (\citeyear{HPS98}) and Park (\citeyear
{P01}) instead of ${\hat{\lambda}}
_{\mathrm{VRS}}$. Theoretical and numerical properties of these testing
procedures are yet to be developed.

\printaddresses

\end{document}